\documentclass[11pt]{amsart}
\usepackage{amsfonts,latexsym,amsthm,amssymb,graphicx}
\usepackage[all]{xy}

\DeclareFontFamily{OT1}{rsfs}{}
\DeclareFontShape{OT1}{rsfs}{n}{it}{<-> rsfs10}{}
\DeclareMathAlphabet{\mathscr}{OT1}{rsfs}{n}{it}

\CompileMatrices

\newtheorem{theorem}{Theorem}[section]
\newtheorem{lemma}[theorem]{Lemma}
\newtheorem{corol}[theorem]{Corollary}
\newtheorem{prop}[theorem]{Proposition}
\newtheorem{claim}[theorem]{Claim}

{\theoremstyle{definition} \newtheorem{defin}[theorem]{Definition}}
{\theoremstyle{remark} \newtheorem{remark}[theorem]{Remark}
\newtheorem{example}[theorem]{Example}}

\newcommand{\Abb}{{\mathbb{A}}}
\newcommand{\Cbb}{{\mathbb{C}}}

\newcommand{\Pbb}{{\mathbb{P}}}
\newcommand{\Zbb}{{\mathbb{Z}}}

\newcommand{\TPbb}{{\Til \Pbb}}

\newcommand{\PGL}{\text{\rm PGL}}

\newcommand{\nG}{\overline{\Pbb}}
\newcommand{\oD}{\overline{D}}
\newcommand{\oE}{\overline{E}}
\newcommand{\oalpha}{{\overline{\alpha}}}
\newcommand{\obeta}{{\overline{\beta}}}
\newcommand{\nor}{n}
\newcommand{\onu}{\overline{n}}

\newcommand{\cC}{{\mathscr C}}
\newcommand{\cD}{{\mathscr D}}
\newcommand{\cCred}{{\cC'}}
\newcommand{\OC}{{\mathscr O_\cC}}

\newcommand{\cS}{{\mathscr S}}
\newcommand{\cX}{{\mathscr X}}
\newcommand{\cY}{{\mathscr Y}}

\newcommand{\Til}[1]{{\widetilde{#1}}}

\newcommand{\Stab}[1]{{\text{Stab}({#1})}}

\DeclareMathOperator{\im}{im}

\title[Limits of PGL(3)-translates of plane curves, II]
{Limits of PGL(3)-translates of plane curves, II}
\author{Paolo Aluffi, Carel Faber}
\address{Dept.~of Mathematics, Florida State University, Tallahassee
FL 32306, U.S.A.}
\email{aluffi@math.fsu.edu}
\address{Inst.~f\"or Matematik, Kungliga Tekniska H\"ogskolan, 
S-100 44 Stockholm, Sweden}
\email{faber@math.kth.se}

\begin{document}

\begin{abstract}
Every complex plane curve $\cC$ determines a subscheme $\cS$ of
the $\Pbb^8$ of $3\times3$ matrices, whose {\em projective normal cone\/} 
(PNC) captures subtle invariants of~$\cC$.

In \cite{ghizzI} we obtain a set-theoretic description of the PNC and
thereby we determine all possible limits of families of plane curves
whose general element is isomorphic to~$\cC$. 
The main result of this article is the determination of the PNC as a 
{\em cycle\/}; this is an essential ingredient in our computation in
\cite{MR2001h:14068} of the degree of the $\PGL(3)$-orbit closure of an
arbitrary plane curve, an invariant of natural enumerative
significance.
\end{abstract}

\maketitle

\section{Introduction}\label{intro}

\subsection{}
In \cite{ghizzI} we determine the possible {\em limits\/} of a fixed, 
arbitrary complex plane curve~$\cC$ obtained by applying to it a family 
of translations $\alpha(t)$ centered at a singular transformation of the 
plane. In other words, we describe the curves in the boundary of the
$\PGL(3)$-orbit closure of a given curve $\cC$.

The germ $\alpha(t)$ is seen as a $\Cbb[[t]]$-valued point of the
projective space $\Pbb^8$ parametrizing $3\times 3$ matrices (up to 
scalar); the determination of the limits of a curve $\cC$ of degree~$d$,
as points of the projective space $\Pbb^N$ of degree-$d$ plane curves,
is accomplished by studying a specific subset of $\Pbb^8\times \Pbb^N$, 
determined by~$\cC$: namely, the projective normal cone (PNC) of the 
subscheme $\cS$ of indeterminacies of the rational map 
$c: \Pbb^8 \dashrightarrow \Pbb^N$ 
associating to $\varphi\in \PGL(3)$ the translate of $\cC$ by 
$\varphi$.

In \cite{ghizzI} we find that the PNC has five types of 
components, reflecting different features of (the support~$\cCred$ of) 
$\cC$: linear components of $\cCred$ (type I); nonlinear components
(type II); singular points at which the tangent cone consists of at least
three lines (type III); sides of the Newton polygon at inflectional or
singular points (type IV); and specially `tuned' formal branches at
singularities of $\cCred$ (type V).
This analysis amounts to a set-theoretic description of the PNC, and suffices 
for the determination of the limits of $\cC$. For applications
to enumerative geometry, and specifically the computation of the degree
of the $\PGL(3)$-orbit closure of $\cC$, it is necessary to have the
more refined information of the PNC {\em as a cycle;\/} that is, the
multiplicities with which the components identified in \cite{ghizzI}
appear in the PNC. This information is listed in \cite{MR2001h:14068},
\S2, and crucially used there as an ingredient in the proof of the
enumerative results. In this article we prove that the multiplicities are as 
stated in \cite{MR2001h:14068}, thereby completing the proof of the results 
in~loc.~cit.

A full statement of the main result of this paper is given in \S\ref{stat};
this section 
also includes a summary of the set-theoretic description given in \cite{ghizzI}.
The general shape of the result is as follows: for each feature of $\cC$ 
`responsible' for a component of the PNC, we give a corresponding
contribution to the multiplicity; the multiplicity for a given component is
obtained by adding these contributions. For example, components of
type IV correspond 
to sides of the Newton polygon for $\cC$ at singular or inflection 
points of its support. We find that the contribution due to a side
of the Newton polygon with vertices $(j_0,k_0)$, $(j_1,k_1)$ (where 
$j_0<j_1$), and corresponding limit for $\cC$
$$x^{\overline q} y^r z^q \prod_{j=1}^S\left(y^c+\rho_j
x^{c-b} z^b\right)\quad,$$
(with $b$ and $c$ relatively prime) equals
$$\frac{j_1k_0-j_0k_1}SA\quad,$$
where $A$ is the number of automorphisms $\Abb^1 \to \Abb^1$, $\rho
\mapsto u\rho$ (with $u$ a root of unity) preserving the $S$-tuple
$\{\rho_1,\dots, \rho_S\}$.
In general, the multiplicities computed here capture delicate information
about the singularities of $\cC$; cf.~especially the multiplicities of components
of type~V. It would be interesting to characterize deformations of $\cC$
for which these multiplicities remain constant.

Explicit examples of computations of multiplicities, for the two curves 
described in Examples~2.1 and~2.2 of \cite{ghizzI}, \S2.6, are given 
in~\S\ref{exas}.

Concerning the proofs: the PNC may be realized as the exceptional divisor
in the blow-up $\TPbb^8$ of $\Pbb^8$ along the subscheme $\cS$ 
mentioned above. The multiplicities of components of type~I and~II,
which depend on `global' features of the curve (the type and multiplicity
of components of the curve) may be computed by analyzing explicitly
parts of this blow-up. Components of type~III, IV, and~V depend
on the behavior of $\cC$ at singularities or inflection points of its support,
and require a more refined analysis. 
This is performed by considering
the normalization $\nG$ of $\TPbb^8$, and studying the pull-back of 
the PNC to $\nG$.
`Marker' germs obtained in \cite{ghizzI} and centered at such 
components may be lifted to germs intersecting the components of this 
pull-back {\em transversally.\/} We use this fact to relate the multiplicity
of a given component $D$ to a weight associated to the corresponding 
marker germ, the number of components $\oD$ of the pull-back 
dominating $D$, and the degrees of the induced maps $\oD \to D$.

The weights may be computed from information collected in \cite{ghizzI}.
The other ingredients are obtained by studying more closely the behavior
of germs $\alpha(t)$ centered at a point of $\cS$, especially in connection 
with the natural $\PGL(3)$-action on~$\TPbb^8$ and~$\nG$. For example, 
every $\alpha(t)$ determines several subgroups of $\PGL(3)$: among 
these are the stabilizer $G$ of the limit $\cX$ of translates of $\cC$ by 
$\alpha(t)$; the stabilizer $\Til G$ of the center $(\alpha, \cX)$ of its lift 
to $\TPbb^8$; and the stabilizer $\overline G$ of the center $\oalpha$ 
of the lift to $\nG$. 
We prove that the degree of $\oD\to D$ equals the index of 
$\overline G$ in $\Til G$ for a corresponding marker germ. 
The group $G$ is available from previous work (cf.~\S1 of 
\cite{MR2002d:14084}); from this it is not hard to obtain a description of
$\Til G$. Further, we identify $\overline G$
with a subgroup of $\PGL(3)$, which we call the `inessential subgroup'
(w.r.t.~$\alpha(t)$),
roughly consisting of elements in the stabilizer of $(\alpha,\cX)$ whose 
presence is due to possible reparametrizations of~$\alpha(t)$.

An explicit computation of the inessential subgroups for the different
types of components then allows us to compute the degree of $\oD
\to D$.

\subsection{Two examples}\label{exas}
We present two examples of computations of multiplicities of all 
components of the PNC for two plane quartic curves. The PNCs for 
these two curves are described (set-theoretically) in 
\cite{ghizzI}, \S2.6, Examples~2.1 and~2.2. The reader should refer to \S\ref{stat}, 
and especially Theorem~\ref{main}, for the notation used in this subsection. 
We also include enumerative consequences for these curves, obtained 
by applying the machinery of \cite{MR2001h:14068} to the information
obtained here.

\begin{example}\label{exone2}
Consider the reducible quartic~$\cC_1$ given by the equation
$$(y+z)(xy^2+xyz+xz^2+y^2z+yz^2)=0\,.$$
It consists of an irreducible nodal cubic and a line through the node
and an inflection point.
In \cite{ghizzI}, Example~2.1, we showed that the PNC for~$\cC_1$ has
seven components:
one each of types~I, II, III, and four components of type~IV.
According to Theorem~\ref{main}, the component of
type I has multiplicity~1 and the component of type II has
multiplicity~2. The multiplicity of the component of type~III
equals $3\cdot6=18$, the product of the multiplicity of the point
and the number of automorphisms of the tangent cone as a triple
in the pencil of lines through the point. To determine the
multiplicities of the four components of type IV, we first note
that in each case $S=1$ and hence $A=1$. For each of the two components
coming from the inflection points, the relevant side of the Newton
polygon has endpoints $(j_0,k_0)=(0,1)$ and $(j_1,k_1)=(3,0)$;
hence their multiplicities equal $j_1k_0-j_0k_1=3$.
For the component coming from the double point and the tangent line
to the cubic,
$(j_0,k_0)=(1,1)$ and $(j_1,k_1)=(4,0)$, so its multiplicity equals~4.
For the component coming from the triple point, both tangent lines
to the cubic give rise to a Newton polygon whose only relevant side has
endpoints $(2,1)$ and $(4,0)$. The multiplicity of this component is
obtained by adding the two contributions; thus it equals $4+4=8$.

The knowledge of the PNC as a cycle makes it possible to apply the
results of \cite{MR2001h:14068},
specifically, Propositions~3.1--3.5. 
One finds that the {\em adjusted
predegree polynomial\/} (a.p.p.)
of~$\cC_1$ equals
$$1+4H+8H^2+\frac{21}{2}H^3+\frac{81}{8}H^4+\frac{147}{20}H^5
+\frac{91}{30}H^6+\frac{89}{140}H^7+\frac{31}{840}H^8.$$
The a.p.p.~is defined on p.~4 of 
\cite{MR2001h:14068}. The `predegree' of the orbit closure of
a curve is $8!\cdot$(the coefficient of $H^8$ in the a.p.p.);
the predegree of~$\cC_1$ equals 
$8!\cdot\frac{31}{840}=1488$.
Since 
the stabilizer of~$\cC_1$ has order~$2$, the degree of its
orbit closure equals $744$.
\qed\end{example}

\begin{example}\label{extwo2}
Consider the irreducible quartic~$\cC_2$ given by the equation
$$(y^2-xz)^2=y^3z.$$
Its singular points are the ramphoid cusp $(1:0:0)$ and the
ordinary cusp $(0:0:1)$; further, it has one inflection point, which
is ordinary. In \cite{ghizzI}, Example~2.2, we showed that the PNC 
for~$\cC_2$ has four components:
one of type~II, two of type~IV, and one of type~V.
According to Theorem~\ref{main}, the component of
type II has multiplicity 2 and the component of type IV coming from the
inflection point has multiplicity 3. For the component of type IV coming
from the ordinary cusp, the relevant side of the Newton polygon has endpoints 
$(0,2)$ and $(3,0)$, so its multiplicity equals~$6$. To compute the
multiplicity of the component of type V coming from the ramphoid cusp,
observe that 
(with notation as in \S\ref{stat})
$C=\frac52$, $f_{(C)}(y)=y^2$, $S=2$, and 
$[\gamma_C^{(1)},\gamma_C^{(2)}]=[1,-1]$, so that
$\ell=1$, $W=\frac52+\frac52=5$, and $A=2\cdot2=4$. Hence the multiplicity
equals $1\cdot5\cdot4=20$. Applying the results of \cite{MR2001h:14068},
one finds that the a.p.p.~of~$\cC_2$ equals
$$1+4H+8H^2+\frac{32}{3}H^3+\frac{32}{3}H^4+\frac{122}{15}H^5
+\frac{61}{16}H^6+\frac{27}{28}H^7+\frac{31}{896}H^8,$$ 
so that its predegree is $1395$
(cf.~Examples~5.2 and~5.4 
of~\cite{MR2001h:14068}: the predegree~$14280$ for the orbit closure of
a general quartic is corrected by~$3960$ because of the ordinary cusp, and
by~$1785\cdot 5=8925$ due to the ramphoid cusp; $14280-3960-8925=1395$).
Since $\cC_2$ has trivial stabilizer, the
degree of its orbit closure equals $1395$ as well.
\qed\end{example}

\subsection{Acknowledgments} Work on this paper was made
possible by support from 
Mathematisches Forschungsinstitut Oberwolfach,
the Volkswagen Stiftung,
the Max-Planck-Institut f\"ur Mathematik (Bonn),
Princeton University,
the G\"oran Gustafsson foundation,
the Swedish Research Council,
the Mittag-Leffler Institute,
MSRI, NSA, NSF, and our home institutions.

\section{Statement of the result}\label{stat}
\subsection{Summary of \cite{ghizzI}}\label{prell}
We begin by recalling briefly the situation examined in~\cite{ghizzI}.

We work over $\Cbb$.
Let $\cC$ be an arbitrary (possibly singular, reducible, non-reduced)
plane curve of degree $d$. The group $\PGL(3)$ acts on
the plane $\Pbb^2$, and hence on the projective space $\Pbb^N=
\Pbb^{d(d+3)/2}$ parametrizing plane curves of degree $d$. Explicitly,
if $F(x,y,z)$ is a homogeneous polynomial generating the ideal of $\cC$,
then $\alpha\in\PGL(3)$ translates $\cC$ to the curve $\cC\circ\alpha$
determined by
$$F(\alpha(x,y,z))\quad.$$
The subset of $\Pbb^N$ consisting of all translates $\cC\circ\alpha$ is 
the {\em linear orbit\/} of $\cC$, which we denote by $\OC$. In \cite{ghizzI} 
we determine the `limits' of $\cC$, by describing the boundary 
($=\overline \OC\smallsetminus\OC$) of the linear orbit of $\cC$.
We do this by analyzing a locus in $\Pbb^8\times \Pbb^N$, components
of which dominate
(most of) the boundary of~$\OC$. 

This locus is the {\em PNC\/} 
determined by $\cC$, as follows. The action map $\PGL(3) \to \Pbb^N$
determined by $\cC$ extends to a rational map
$$c: \Pbb^8 \dashrightarrow \Pbb^N$$
defined by $\alpha \mapsto \cC\circ\alpha$ for $\alpha\in \Pbb^8$ a
$3\times 3$ matrix. We let $\cS$ be the scheme of indeterminacies
of this rational map. The PNC (projective normal cone) of $\cS$ in
$\Pbb^8$ is a purely 7-dimensional scheme $E$, which can be naturally
embedded in $\Pbb^8\times \Pbb^N$.

The set-theoretic description of the 
PNC in \cite{ghizzI} identifies five `types' of components, arising from 
different features of (the support $\cCred$ of) $\cC$:
\begin{itemize}
\item[I:] The linear components;
\item[II:] The nonlinear components;
\item[III:] The singular points at which the tangent cone
is supported on $\ge 3$ lines;
\item[IV:] The Newton polygons at its singular and inflection points;
\item[V:] The Puiseux expansions of formal branches
at its singular points.
\end{itemize}
Each of the corresponding components of the PNC may be described 
by giving general points on it; the component is obtained as the orbit
closure of such a point under the natural right action of $\PGL(3)$ on 
$\Pbb^8\times\Pbb^N$, or as the closure of a union of such orbits
in the type~II case.

These general points $(\alpha,\cX)\in \Pbb^8\times\Pbb^N$ are as listed 
below (cf.~\cite{ghizzI}, \S2.3).
Concerning our terminology, a {\em star\/} of lines is a set of concurrent 
lines; a {\em fan\/} is a star union a general line; 
the reader is addressed to
\cite{MR2002d:14084}, \S1, or the summary in \cite{ghizzI}, \S6,
for a more extensive description of the 
curves appearing in this list.
\begin{itemize}
\item[I:] $\im\alpha=$ a line $L\subset\cC$; $\cX=$ a star or fan centered 
at $\ker\alpha$;
\item[II:] $\im\alpha=$ a general point of a nonlinear component $\cD$
of $\cCred$;
$\cX=$ a (possibly multiple) conic, (possibly) union a multiple tangent line
supported on $\ker\alpha$;
\item[III:] $\im\alpha=$ a singular point of multiplicity $m$ on $\cC$ at which 
the tangent cone consists of three or more lines; $\cX=$ a star union
a line of multiplicity $d-m$ supported on $\ker\alpha$;
\item[IV:] $\im\alpha=$ a singular or inflection point of $\cC$ at which the 
Newton polygon of $\cC$ has a segment of slope strictly between $-1$ 
and $0$; $\cX=$ a union of (possibly degenerate, possibly multiple) 
cuspidal curves;
\item[V:] $\im\alpha=$ a point of $\cC$ at which the formal branches
of $\cCred$ have suitable truncations $f_{(C)}(y)=
\sum_{\lambda_i<C} \gamma_{\lambda_i} y^{\lambda_i}$
depending on certain rational numbers (`characteristics') $C$; $\cX=$ a 
union of quadritangent conics, union (possibly) a multiple tangent line 
supported on $\ker\alpha$.
\end{itemize}

The limits in type~IV and~V may be
written explicitly in the form given below in Theorem~\ref{main}.

\subsection{The main theorem}\label{multstat}
The goal of this paper is to extend the set-theoretic description of
the PNC obtained in \cite{ghizzI} and recalled above 
to a description of the PNC as a {\em cycle,\/} that is, to determine
the multiplicities of the various components.

For each feature of $\cC$ 
responsible for a component of the PNC, we give a corresponding
contribution to the multiplicity; the multiplicity for a given component is
obtained by adding these contributions.

The main result of this paper is the following (compare 
with~\cite{MR2001h:14068}, \S2, Facts~1 through~5).

\begin{theorem}\label{main}
Each feature of $\cC$ contributes a multiplicity to the corresponding 
component, as follows:

\begin{itemize}

\item{\em Type I.\/} The multiplicity of the component determined by a
  line $L\subset \cC$ equals the multiplicity of $L$ in
  $\cC$.

\item{\em Type II.\/} The multiplicity of the component determined by
  a nonlinear component $\cD$ of $\cC$ equals $2m$,
  where $m$ is the multiplicity of $\cD$ in $\cC$.

\item{\em Type III.\/} The multiplicity of the component determined by
  a singular point $p$ of $\cC$ such that the tangent cone
  $\lambda$ to $\cC$ at $p$ is supported on three or more lines
  equals $mA$, where $m$ is the multiplicity of $\cC$ at $p$
  and $A$ equals the number of automorphisms of $\lambda$ as an 
  $m$-tuple
  in the pencil of lines through $p$.

\item{\em Type IV.\/} The multiplicity of the component determined by
  one side of a Newton polygon for $\cC$, with endpoints
  $(j_0,k_0)$, $(j_1,k_1)$ (where $j_0<j_1$) and limit
$$x^{\overline q} y^r z^q \prod_{j=1}^S\left(y^c+\rho_j
x^{c-b} z^b\right)\quad,$$
(with $b$ and $c$ relatively prime) equals
$$\frac{j_1k_0-j_0k_1}SA\quad,$$
where $A$ is the number of automorphisms $\Abb^1 \to \Abb^1$, $\rho
\mapsto u\rho$ (with $u$ a root of unity) preserving the $S$-tuple
$[\rho_1,\dots, \rho_S]$.\footnote{Note: the number $A$ given here
  is denoted $A/\delta$ in \cite{MR2001h:14068}.}

\item{\em Type V.\/} The multiplicity of the component corresponding
  to the choice of a characteristic $C$ and a formal branch 
  $f(y)=\sum_{i\ge 0} \gamma_{\lambda_i} y^{\lambda_i}$
  at a point $p$, with limit
$$x^{d-2S}\prod_{i=1}^S\left(zx-\frac {\lambda_0(\lambda_0-1)}2
\gamma_{\lambda_0}y^2 -\frac{\lambda_0+C}2
\gamma_{\frac{\lambda_0+C}2}yx-\gamma_C^{(i)}x^2\right)$$
is $\ell WA$, where:
\begin{itemize}
\item $\ell$ is the least positive integer $\mu$ such that
  $f_{(C)}(y^\mu)$ has integer exponents.
\item $W$ is defined as follows. For each formal branch $\beta$ of
  $\cC$ at $p$, let $v_\beta$ be the first exponent at which
  $\beta$ and $f_{(C)}(y)=
\sum_{\lambda_i<C} \gamma_{\lambda_i} y^{\lambda_i}$
differ, and let $w_\beta$ be the minimum of
  $C$ and $v_\beta$. Then $W$ is the sum $\sum w_\beta$.
\item $A$ is twice the number of automorphisms $\gamma\to u\gamma+v$
preserving the $S$-tuple $[\gamma_C^{(1)},\dots,\gamma_C^{(S)}]$.
\end{itemize}

\end{itemize}
\end{theorem}

Note that a given component may be obtained in different ways
(each producing a multiplicity, given in the statement).
Concerning components of type I, II, or III, there is
only one contribution for each of the specified data---that is,
exactly one contribution of type~I from each line contained in
$\cC$, one contribution of type~II from each nonlinear
component $\cD$ of $\cC$, and one of type~III from each singular
point of $\cC$ at which the tangent cone is supported on three
or more distinct lines.

As usual, the situation is a little more complex for components of
type~IV and~V. The following information is necessary in order to
apply Theorem~\ref{main} to a given curve~$\cC$.

{\em Components of type~IV\/} correspond to sides of Newton polygons;
one polygon is obtained for each line in the tangent cone at a fixed
singular or inflection point $p$ of $\cC$, and each of these polygons
provides a set of sides (with slope strictly between $-1$ and~$0$).
Exactly one contribution has to be counted for each side obtained in
this fashion. Note that sides of different Newton polygons (for different
lines in the tangent cone of $\cC$ at the same point) may lead to
the same limits by germs centered at the same point, hence to the 
same component of the PNC.

{\em Components of type~V\/} are determined by a choice of a singular
point $p$ of $\cC$, a line~$L$ in the tangent cone to $\cC$ at $p$, a 
characteristic $C$ and a formal branch $z=f(y)$ of $\cC$,
tangent to $L$. Recall (from \cite{ghizzI}, \S2.3 and~2.4)
that these data determine a triple of positive
integers $a<b<c$ with $C=c/a$. Again, different choices may
lead to the same component of the PNC, and we have to specify when
choices should be counted as giving separate contributions. Of course
different points $p$ or different lines in the tangent cone at $p$
give separate contributions; the question is when two sets of data
$(C,f(y))$ for the same point, with respect to the same tangent
line, should be counted separately.

We say that $(C,f(y))$, $(C',g(y))$ are {\em sibling\/} data 
if $C=C'$, the corresponding triples of positive integers are identical,
and the truncations $f_{(C)}(t^a)$, $g_{(C)}(t^a)$ are related by
$g_{(C)}(t^a)=f_{(C)}((\xi t)^a)$ for an $a$-th root $\xi$ of $1$.
(Note that $f_{(C)}((\xi t)^a)$ does not equal $f_{(C)}(t^a)$ in general, since
formal branches may have fractional exponents.)

Then: two pairs $(C,f(y))$, $(C',g(y))$ at the same point, 
with respect to the same tangent line, yield separate contributions if and only 
if they are {\em not\/} siblings.
\vskip 6pt

The proof of Theorem~\ref{main} occupies the rest of 
this paper. Rather direct arguments can be given for the components 
of type~I and type~II, depending on `global' features of~$\cC$; we treat 
these cases in \S\ref{IandII}. The `local' types III, IV, and~V require the
development of appropriate tools, and are treated in \S\S\ref{local}--\ref{Eop}.

\section{Components of type I and II}\label{IandII}
\subsection{Type~I}\label{multtypeI}
The PNC of $\cC$ is denoted by $E$; recall that it may be
identified with the exceptional divisor of the blow-up of $\Pbb^8$ along
the subscheme $\cS$ (cf.~\S\ref{prell}). 

\begin{prop}[Type~I]\label{typeImul}
Assume that $\cC$ contains a line $L$ with multiplicity $m$, and let 
$(\alpha, \mathcal X)$ be a general point of the corresponding component 
$D$ of $E$. Then $\TPbb^8$ is nonsingular at $(\alpha,\mathcal X)$, 
and $D$ appears with multiplicity $m$ in $E$. 
\end{prop}

\begin{proof}
We are going to show that, in a neighborhood of $(\alpha, \mathcal
X)$, $\TPbb^8$ is isomorphic to the blow-up of $\Pbb^8$ along
the $\Pbb^5$ of matrices whose image is contained in $L$. 
The nonsingularity of $\TPbb^8$ near $(\alpha, \mathcal X)$ follows
from this.

Choose coordinates so that $L$ is the line $z=0$, and
$$\alpha=\begin{pmatrix}
1 & 0 & 0 \\
0 & 1 & 0 \\
0 & 0 & 0
\end{pmatrix}\quad;$$
consider the open neighborhood $U\subset\Pbb^8$ of $\alpha$, with 
coordinates
$$\begin{pmatrix}
1 & p_1 & p_2 \\
p_3 & p_4 & p_5 \\
p_6 & p_7 & p_8
\end{pmatrix}\quad.$$
The $\Pbb^5$ of matrices with image contained in $L$ intersects this open 
set along $p_6=p_7=p_8=0$; hence we can choose coordinates 
$q_1,\dots,q_8$ in an affine open subset $V$ of the blow-up of $\Pbb^8$ 
along $\Pbb^5$ so that the blow-up map is given by
$$\left\{
\aligned
p_i &=q_i\qquad i=1,\dots,5 \\
p_6 &=q_6 \\
p_7 &=q_6 q_7\\
p_8 &=q_6 q_8
\endaligned\right.$$
(the part of the blow-up over $U$ is covered by three such open sets;
it will be clear from the argument that the choice made here is
immaterial). With these coordinates, the exceptional divisor has equation 
$q_6=0$.

Under the hypotheses of the statement, the ideal of $\cC$ is
generated by $z^m G(x,y,z)$, where $z$ does not divide $G$; that is,
$G(x,y,0)\not\equiv 0$. 
The rational map $c:\Pbb^8 \dashrightarrow \Pbb^N$ sends 
$(p_1,\dots,p_8)\in U$ to the curve with ideal generated by
$$(p_6 x+p_7 y+p_8 z)^m G(x+p_1 y+p_2 z,p_3 x+p_4 y+p_5 z, p_6 x+p_7
y+p_8 z)\quad;$$
it follows that the ideal of $\cS$ in $U$ is generated by these polynomials
(in $p_1,\dots,p_8$) as $(x:y:z)$ varies in $\Pbb^2$. 
Composing with the blow-up map:
$$V \longrightarrow U \dashrightarrow \Pbb^N\quad,$$
these generators pull-back to
$$q_6^m (x+q_7 y+q_8 z)^m G(x+q_1 y+q_2 z,q_3 x+q_4 y+q_5 z, 
q_6 x+q_6 q_7 y+q_6 q_8 z)\quad.$$
By the hypothesis on $G$, this shows that, along a dense open set $W$
of $V$ intersecting the exceptional divisor $q_6=0$, the ideal of $\cS$ 
pulls back to $(q_6^m)$;
 by the universal property of blow-ups we obtain an induced map
$$W \to \TPbb^8$$
mapping the exceptional divisor $q_6=0$ to $D$. A coordinate
verification shows that this is an isomorphism onto the
image in a neighborhood of a general point of the exceptional divisor,
proving that $\TPbb^8$ is nonsingular in a neighborhood of a
general $(\alpha, \mathcal X)$ in $D$, and that the multiplicity of
$D$ in $E$ is $m$, as stated.
\end{proof}

Proposition~\ref{typeImul} yields the multiplicity statement
concerning type~I components in Theorem~\ref{main}; also cf.~Fact~2~(i)
in \S2 of \cite{MR2001h:14068}.

\subsection{Type~II}\label{multtypeII} 
Components of type~II correspond to nonlinear components $\cD$
of $\cCred$.

\begin{prop}[Type~II]\label{typeIImul}
Assume a nonlinear component $\cD$ appears with multiplicity $m$
in $\cC$, and let $D$ be the corresponding component of $E$.
Then $D$ appears with multiplicity $2m$ in $E$.
\end{prop}

This can be proved by using the blow-ups described in
\cite{MR94e:14032}, which resolve the indeterminacies of the basic
rational map $\Pbb^8 \dashrightarrow \Pbb^N$ over nonsingular, 
non-inflectional points
of $\cC$. We sketch the argument here, leaving detailed
verifications to the reader.

\begin{proof}
In \cite{MR94e:14032} it is shown that two blow-ups at smooth centers
suffice over nonsingular, non-inflectional points of $\cC$. While the curve 
was assumed to be reduced and irreducible in loc.~cit., the reader may 
check that the same blow-ups resolve the indeterminacies over a possibly 
multiple component $\cD$, near nonsingular, non-inflectional points of the 
support of $\cD$. Let $V$ be the variety obtained after these two blow-ups.

Since the basic rational map is resolved by $V$ over a general point
of $\cD$, the inverse image of the base scheme $\cS$ is locally
principal in $V$ over such points. By the universal property of
blow-ups, the map $V \to \Pbb^8$ factors through $\TPbb^8$ over a
neighborhood of a general point of $\cD$.
It may then be checked that the second exceptional divisor obtained in
the sequence maps birationally onto~$D$, and appears with a
multiplicity of $2m$.
The statement follows.
\end{proof}

Proposition~\ref{typeIImul} yields the multiplicity statement
concerning type~II components in Theorem~\ref{main}; also 
cf.~Fact~2~(ii) in \S2 of \cite{MR2001h:14068}.

\section{Components of type III, IV, and V}\label{local}
\subsection{Normalizing the graph}\label{beginprel}
The computation of the multiplicity of components of type III, IV, and V 
is considerably subtler, and requires some preparatory work.

Our main tool will be the {\em normalization\/}
$\nG$ of the closure $\TPbb^8$ of the graph of the basic rational
map $c:\Pbb^8 \dashrightarrow \Pbb^N$ from \S\ref{stat}. We
denote by
$$\nor: \nG \rightarrow \TPbb^8$$ 
the normalization map, and by $\onu$ the composition $\nG \to
\TPbb^8 \to \Pbb^8$.

Recall that the PNC may be realized as the exceptional divisor 
$E$ in the blow-up $\TPbb^8$ of $\Pbb^8$ along the scheme $\cS$ of
indeterminacies of $c$. If $F\in\Cbb[x,y,z]$ generates the ideal of
$\cC$ in $\Pbb^2$, then the ideal of $\cS$ in $\Pbb^8$ is
generated by all expressions
$$F(\varphi(x_0,y_0,z_0))\quad,$$
viewed as polynomials in $\varphi\in\Pbb^8$, as $(x_0,y_0,z_0)$
ranges over $\Pbb^2$. The ideals of $E$ and of
$\oE=\nor^{-1}(E)=\onu^{-1}(\cS)$ are generated by the pull-backs of 
$F(\varphi(x_0,y_0,z_0))$ to $\TPbb^8$, respectively~$\nG$.
 
Denote by $E_i$ the supports 
of the components of $E$, and by  $m_i$ the multiplicity of $E_i$ in $E$.
Also, denote by 
$$\oE_{i1}\,,\,\dots\,,\,\oE_{ir_i}$$
the supports of the components of $\oE$ lying above a given component
$E_i$ of $E$. Finally, let $m_{ij}$ be the multiplicity of
$\oE_{ij}$ in $\oE$. Summarizing: as cycles,
$$[E]=\sum m_i [E_i]\quad,\quad [\oE]=\sum m_{ij} [\oE_{ij}]\quad.$$

\begin{lemma}\label{multcount}
We have
$$m_i=\sum_{j=1}^{r_i} e_{ij} m_{ij}$$
where $e_{ij}$ is the degree of the map $\nor|_{\oE_{ij}}:\oE_{ij} \to E_i$.
\end{lemma}

\begin{proof}
Use the projection formula (\cite{85k:14004}, 
Prop.~2.3~(c)) 
and $(\nor|_{\oE_{ij}})_* [\oE_{ij}]=e_{ij}[E_i]$.
\end{proof}

\begin{example}
There is exactly one component of $\oE$ for each component of type
I~or~II of $E$, mapping birationally to such a component. 
This may be established by analyzing the blow-ups in 
Propositions~\ref{typeImul} and~\ref{typeIImul}.
\end{example}

\subsection{Three propositions}\label{result}
The rest of the paper consists of the computation of the ingredients
needed in order to apply Lemma~\ref{multcount} to the case of
components of type III, IV, and V. It is not hard to
extract the multiplicities $m_{ij}$ from more refined information
collected in \cite{ghizzI}.
The number of components of $\oE$ dominating 
a given component of $E$, and the degrees $e_{ij}$, will require more work.

The following three propositions collect the results we will obtain. Proofs
of these results are presented in \S\ref{mulnorm} and ff.

Concerning components of type~III, the situation is very simple:
\begin{prop}[Type~III]\label{typeIIImul}
Let $D$ be a component of type~III, corresponding to a point $p$
of multiplicity $m$
at which the tangent cone to $\cC$ is supported on at least $3$ lines.
\begin{itemize}
\item There is exactly one component $\oD$ of $\oE$ dominating $D$.
\item The degree of the map $\oD \to D$ equals the number of linear
automorphisms of the $m$-tuple determined by the tangent cone to $\cC$
at $p$.
\item The multiplicity of $\oD$ equals $m$.
\end{itemize}
\end{prop}

As recalled in \S\ref{prell}, components of type~IV correspond to the
choice of a point $p\in \cC$, a line $L$ in the tangent cone to $\cC$
at $p$, and one side of the Newton polygon for $\cC$ at $p$ relative
to $L$, of slope $-b/c$ with $0<b<c$.
The same component may arise from sides with the same slope,
with respect to different lines in the tangent cone at~$p$. Limits
corresponding to these choices are of the form
\begin{equation}\tag{*}
x^{\overline q}y^rz^q \prod_{j=1}^S(y^c+\rho_j x^{c-b}z^b)\quad,
\end{equation}
with $\rho_j \ne 0$,
where $S+1$ is the number of lattice points on the chosen side of the
Newton polygon.

\begin{prop}[Type~IV]\label{typeIVmul}
Let $D$ be a component of $E$ of type~IV, as above. Then:
\begin{itemize}
\item There is exactly one component $\oD$ of $\oE$ over~$D$ for 
each line $L$ in the tangent cone to $\cC$ at $p$, with respect to 
which the Newton polygon has a side of slope $-b/c$, that leads to limit (*).

\item The degree of the map $\oD \to D$ equals the number of 
automorphisms $\Abb^1 \to \Abb^1$, 
$\rho\mapsto u\rho$ (with $u$ a root of unity)
 preserving the 
$S$-tuple $[\rho_1,\dots,\rho_S]$.

\item Let $(j_0,k_0)$ and $(j_1,k_1)$
be the endpoints of the side of the Newton polygon $(j_0<j_1)$.
The multiplicity of $\oD$ is
$$\frac{j_1k_0-j_0k_1}S\quad.$$
\end{itemize}
\end{prop}

Components of type~V are determined as follows. Choose a point $p\in\cC$,
a line $L$ in the tangent cone to $\cC$ at $p$, and coordinates so that
$p=(1:0:0)$, $L$ is the line with equation $z=0$, and $y=0$ is not part
of the tangent cone to $\cC$ at $p$.
Express $\cC$ at $p$ in terms of {\em formal branches\/} 
(cf.~\S4.1 in \cite{ghizzI}),
with Puiseux expansions of the form
$$z=f(y) = \sum_{i\ge 0} \gamma_{\lambda_i}y^{\lambda_i}\quad,$$
where $\lambda_0\ge 1$, and $\lambda_i<\lambda_{i+1}$.

These choices determine a finite set of positive rational numbers (called
`characteristics' in \cite{ghizzI}):
that is, those numbers $C$ which are exponents $\lambda_i$, $i>0$, for some 
formal branch tangent to $L$; and such that at least two formal branches
have the same truncation `modulo $y^C$':
$$f_{(C)}(y) = \sum_{\lambda_i<C} \gamma_{\lambda_i}y^{\lambda_i}
\quad,$$
but different coefficients for $y^C$. Note that $C>\lambda_0$.

For fixed $p$ and $L$,
the choice of a characteristic $C$ and of one such
formal branch
determines a type~V component.

If $z=f(y) = \sum_{i\ge 0} \gamma_{\lambda_i}y^{\lambda_i}$ is a coordinate
representation of the formal branch (with $\lambda_0>1$), the limit 
corresponding to the choice 
of $p$, $L$, and $(C,f(y))$ is (cf.~\S\ref{multstat})
\begin{equation}\tag{**}
x^{d-2S}\prod_{i=1}^S\left(zx-\frac {\lambda_0(\lambda_0-1)}2
\gamma_{\lambda_0}y^2 -\frac{\lambda_0+C}2
\gamma_{\frac{\lambda_0+C}2}yx-\gamma_C^{(i)}x^2\right)\quad,
\end{equation}
where $\gamma_C^{(i)}$ are the coefficients of $y^C$ for all $S$ formal
branches sharing the truncation.

We note that, for a fixed point $p$, different lines $L$ and different $(C,f(y))$
may lead to the same limit, and hence to the same component $D$ of $E$.
While it is clear that different lines must correspond to different components
$\oD$ of $\oE$ over $D$, the question of which pairs $(C,f(y))$ correspond
to the same component $\oD$ is subtle, and accounts for the most technical
parts of this paper.

The choice $(C,f(y))$ determines three integers
$a<b<c\,$: let $B=\frac{C-\lambda_0}2+1$ (so that $1<B<C$, as 
$C>\lambda_0> 1$), and let
$$a=\text{least positive integer such that $aB$, $aC$, and all
$a\lambda_i$ for $\lambda_i<C$ are integers}\quad;$$
we then set $b=aB$, $c=aC$.

We say that $(C,f(y))$ as above and $(C',g(y))$ are {\em sibling\/} 
data if  the corresponding integers $a<b<c$, $a'<b'<c'$ are the same (so in
particular $C=C'$) and further 
$$g_{(C)}(y)=\sum_{\lambda_i<C} \xi^{a\lambda_i} \gamma_{\lambda_i}
y^{\lambda_i}$$
for an $a$-th root $\xi$ of $1$. The right-hand side is well-defined
because $a\lambda_i\in \Zbb$ for $\lambda_i<C$.
Siblings lead to the same component of $\oE$ (cf.~Claim~\ref{siblingsame}).
\vskip 6pt

{\em Remark.}
If $z=f(y)$, $z=g(y)$ are formal branches belonging to
the same irreducible branch of $\cC$ at $p$, then the
corresponding data $(C,f(y))$, $(C,g(y))$ 
are siblings for all $C$.
Indeed, if the branch has multiplicity $m$ at $p$ then
$f(\tau^m)=\varphi(\tau)$ and $g(\tau^m)=\psi(\tau)$, with
$\psi(\tau)=\varphi(\zeta\tau)$
for an $m$-th root $\zeta$ of $1$ (\cite{MR1836037}, \S7.10). That is,
$$\text{if}\quad 
f(y) =\sum \gamma_{\lambda_i} y^{\lambda_i}\quad,\quad\text{then}\quad
g(y) =\sum \zeta^{m\lambda_i}\gamma_{\lambda_i} y^{\lambda_i}$$
for $\zeta$ an $m$-th root of $1$. 
Now let $\rho$ be an $(am)$-th root of $1$ such that $\rho^a=\zeta$,
and set $\xi=\rho^m$; since the exponents $a\lambda_i$ in the
truncations are integers, and so are all exponents 
$m\lambda_i$, we have
$$\zeta^{m\lambda_i}=\rho^{ma\lambda_i}=\xi^{a\lambda_i}$$
for all exponents $\lambda_i<C$, and this shows that the truncations
are siblings.
\vskip 6pt

For a given $p$ and $L$, the set of $(C,f(y))$ leading to a given
component is partitioned into sibling classes. By the remark, 
all formal branches belonging to a given irreducible branch of
$\cC$ tangent to $L$ at $p$ are in the same class.
The sibling classes can therefore be thought of as particular
collections of irreducible branches with a common tangent.

We let $A$ be the number of components of the stabilizer of the limit (**);
that is, by \cite{MR2002d:14083}, \S4.1, twice the number of
automorphisms $\gamma\mapsto u\gamma+v$ preserving the $S$-tuple
$[\gamma_C^{(1)},\dots,\gamma_C^{(S)}]$.

Further, we let $h$ denote the greatest common divisor of $a$ and all
$a\lambda_i$ for $\lambda_i<C$.

Finally, for every choice of $L$, $C$, and $f(y)$, and every
formal branch $\beta$ of $\cC$ at $p$, define a rational number 
$w_\beta$ as follows:
\begin{itemize}
\item if the branch is not tangent to $L$, then $w_\beta=1$;
\item if the branch is tangent to the line $L$, but does not
  truncate to $f_{(C)}(y)$, then $w_\beta=$ the first exponent at which
  $\beta$ and $f_{(C)}(y)$ differ;
\item if the branch truncates to $f_{(C)}(y)$, then $w_\beta=C$.
\end{itemize}
Note that $aw_\beta$ is an integer for all $\beta$. 
We let $W$ denote the sum $\sum w_\beta$.

\begin{prop}[Type~V]\label{typeVmul}
Let $D$ be a component of $E$ of type~V, determined by the choice
of $p$ and a limit (**). Then:
\begin{itemize}
\item The set of components $\oD$ of $\oE$ over~$D$ is in bijection with
the set of all sibling classes contributing to $D$ (for all lines~$L$).
\item For a choice of a line $L$ and of $(C,f(y))$, the degree of the
map $\oD \to D$ for the corresponding component $\oD$ equals $\frac Ah$.
\item The multiplicity of $\oD$ equals $aW$.
\end{itemize}
\end{prop}

The statements about multiplicities
in Propositions~\ref{typeIIImul}---\ref{typeVmul} may be summarized 
as follows. Components $\oD$ of $\oE$ will correspond to germs in
a standard form, to be introduced in \S\ref{markergerms}. For a general 
$q\in \Pbb^2$, consider the (parametrized) curve $\cY$ obtained by 
applying to $q$ one of these germs. Then the multiplicity of $\oD$ 
in $\oE$ is the intersection multiplicity of $\cY$ and $\cC$ at $p$.

\subsection{Propositions~\ref{typeIIImul}---\ref{typeVmul} imply 
the main theorem}
Propositions~\ref{typeIIImul} and~\ref{typeIVmul} imply the statements for 
type III and IV in Theorem~\ref{main}, as an immediate consequence
of Lemma~\ref{multcount}.

It is perhaps less evident that Proposition~\ref{typeVmul} implies
the formula for the multiplicity of a type~V component given in 
Theorem~\ref{main}: according to Proposition~\ref{typeVmul}
and Lemma~\ref{multcount}, the sibling class of $(C,f(y))$
gives a contribution of $a W\frac Ah$ to the multiplicity of the
corresponding component of type~V; we have to check that
we have $\ell=\frac ah$.

For this, let $\lambda_i$, $i=1,\dots,r$ be the exponents appearing in
$f_{(C)}(y)$. If $h'$ is any divisor of $a$ and all $a\lambda_i$, then
as $\frac a{h'}\lambda_i$ are integers, necessarily $\frac a{h'}$ is a
multiple of $\ell$. That is, $h'$ divides $\frac a\ell$. On the other
hand, $\frac a\ell$ is a divisor of $a$ and all $a\lambda_i$. Hence
$\frac a\ell$ equals the greatest common divisor of $a$ and all
$a\lambda_i$, that is, $h$, as needed.

Summarizing, we are reduced to proving
Propositions~\ref{typeIIImul}, ~\ref{typeIVmul}, and~\ref{typeVmul}. 
The reader should compare the statements of these propositions 
with~\cite{MR2001h:14068}, \S2, Facts~3 through~5.
\smallskip

The proof of Propositions~\ref{typeIIImul}, ~\ref{typeIVmul}, 
and~\ref{typeVmul} is organized as follows. The multiplicities of the
components of $\oE$ are determined in \S\ref{mulnorm}, in terms
of the weights of the marker germs found in \cite{ghizzI}.
In \S\ref{compnor} we enumerate the components of $\oE$ over a 
given component of the PNC; the main tool is obtained in
Corollary~\ref{critesamecom}, and its application to our situation 
(especially for components of type~V) relies on the technical
Lemma~\ref{technical}.
Finally, in \S\ref{Eop} we compute the degree of a component of 
$\oE$ over the corresponding component of the PNC. The main tool
here is Proposition~\ref{degreetool}, which relates this degree to the
`inessential subgroup' mentioned in \S\ref{intro}. The inessential
subgroups for components of type~III, IV, and~V are determined 
in~\S\ref{degoDoD}, concluding the proof.

\section{Marker germs, and multiplicities in the normalization}\label{mulnorm}
\subsection{Marker germs}\label{markergerms}
The statements about multiplicities of components in the normalization
will be straightforward consequences of more refined information obtained in
\cite{ghizzI}; we begin by recalling this information.

Every germ $\alpha(t)$ in $\Pbb^8$, whose general element is
invertible, and such that $\alpha(0)\in \cS$, lifts to a unique germ in
$\TPbb^8$ centered at a point $(\alpha(0),\cX)$ of the PNC. 
The germ lifts to a unique germ $\oalpha(t)$ in $\nG$, centered at a 
point $\oalpha=\oalpha(0)$ of $\oE$. Conversely, every germ in $\nG$ that is 
not contained in one of the components $\oE_{ij}$ is the lift of a unique germ 
in $\Pbb^8$.

The data obtained in \cite{ghizzI} includes a list of {\em marker\/}
germs, marking components of different types. For types III, IV, 
and V these are as follows (cf.~\S2.3 in \cite{ghizzI}). 

Components of type~III, IV, and~V determine 
a point $p$ of $\cC$; choose coordinates so that this point is 
$(1:0:0)$. Type~IV and~V depend on the choice of a line $L$ in
the tangent cone to $\cC$ at $p$; choose coordinates so that
this line is the line~$z=0$.

{\bf Type III.\/}
The corresponding marker germ is
$$\alpha_{III}(t)=\begin{pmatrix}
1 & 0 & 0 \\
0 & t & 0 \\
0 & 0 & t
\end{pmatrix}\quad.$$

{\bf Type IV.\/}
These components are determined by the choice of a side with slope
strictly between $-1$ and $0$ of the Newton 
polygon for $\cC$ at $p$, with respect to $L$. Let $b$ and $c$ be relatively 
prime positive integers, such that  $-b/c$ is this slope. Then the corresponding
marker germ is
$$\alpha_{IV}(t)=\begin{pmatrix}
1 & 0 & 0 \\
0 & t^b & 0 \\
0 & 0 & t^c
\end{pmatrix}\quad.$$

{\bf Type V.\/}
These components are determined by the choice of a formal branch
$z=f(y)=\gamma_{\lambda_0}y^{\lambda_0}+\dots$ for $\cC$ at $p$ 
tangent to $L$, and of a characteristic $C>\lambda_0$.
For $a<b<c$ positive integers such that $\frac ca=C$, $\frac ba=
\frac{C-\lambda_0}2+1$, the corresponding marker germ is
$$\alpha_V(t)=\begin{pmatrix}
1 & 0 & 0 \\
t^a & t^b & 0 \\
\underline{f(t^a)} & \underline{f'(t^a)t^b} & t^c
\end{pmatrix}$$
where $\underline{\cdots}$ denotes the truncation modulo $t^c$. 
The integer $a$ is chosen to be the minimum one for which all entries in this 
germ are polynomials.

\subsection{Equivalence of germs}
In \cite{ghizzI}, \S3.1.1, we consider the following notion of `equivalence' of germs:

\begin{defin}\label{equivgermsnew}
Two germs $\alpha(t)$, $\beta(t)$ are {\em equivalent\/} if 
$\beta(t\nu(t))\equiv \alpha(t)\circ m(t)$, with $\nu(t)$ a unit
in $\Cbb[[t]]$, and $m(t)$ a germ such that $m(0)=I$.
\end{defin}

Intuitively speaking, equivalent germs may be `deformed continuously'
one into the other, while keeping their center fixed. This notion will be 
crucial in the rest of the paper. We will first prove that two germs lift 
to $\nG$ 
to germs with the
same center if they are equivalent; in essence, the converse also 
holds (cf.~Proposition~\ref{criteqlift}).

\begin{lemma}\label{equivgermso}
If $\alpha(t)$, $\beta(t)$ are equivalent germs, then $\oalpha=\obeta$.
\end{lemma}

\begin{proof}
Since the center of the lift does not depend on a change of parameter,
we may assume $\beta(t)\equiv \alpha(t)\circ m(t)$, with $m(t)$ a germ
centered at the identity $I$. Let $F(x,y,z)$ be a generator of the
ideal of $\cC$, and let
$$\alpha_h(t)=\alpha(t)\circ((1-h)\, I+h\, m(t))\quad.$$
Then $\alpha_h(t)$ is equivalent to $\alpha(t)$ for all $h$; in particular,
the initial term of $F\circ\alpha_h(t)$ is independent of $h$
(cf.~Lemma~3.2 in \cite{ghizzI}),
so the center of the lift of $\alpha_h(t)$ to $\TPbb^8$ is independent of 
$h$, and it follows that $\oalpha_h\in \nG$ is independent of $h$ as
$n:\nG \to \TPbb^8$ is a finite map.
The statement follows, since $\oalpha=\oalpha_0$ and 
$\obeta=\oalpha_1$.
\end{proof}

A second reason why the notion of equivalence is important in this paper,
as well as in \cite{ghizzI}, is the following fact.
A germ `contributes' a component to $E$ if its lift to~$\TPbb^8$ is
centered at a general point of that component (a more
precise definition will be given in~\S\ref{pgl3act}).

\begin{lemma}\label{markequiv}
Every contributing germ $\alpha(t)$ is equivalent to a marker germ
(in suitable coordinates, and possibly up to replacing $t$ by a root
$t^{1/k}$).\qed
\end{lemma}

\begin{remark}\label{coordchoice}
The coordinate choices implicit in this statement are important, and
we discuss them here. 
Lemma~\ref{markequiv} is proved (cf.~\cite{ghizzI}, Lemma~3.6 and ff.) by 
first showing that every germ $\alpha(t)$ may be written as
$$\alpha(t)=H\cdot \begin{pmatrix}
1 & 0 & 0 \\
q & 1 & 0 \\
r & s & 1 \end{pmatrix}\cdot \begin{pmatrix}
1 & 0 & 0 \\
0 & t^b & 0 \\
0 & 0 & t^c
\end{pmatrix}
\cdot m(t)\quad,$$
where $m(t)$ is invertible, and $q$, $r$, $s$ are polynomials satisfying 
certain conditions;
for example, $q(t)$ may be assumed to be either $0$, or a power $t^a$
(possibly after a parameter change).
Coordinates are then chosen in the plane containing $\cC$
 so that $H=I$, and it is shown that
the hypothesis that $\alpha(t)$ contributes a component forces certain 
conditions on $b$, $c$, and $q$, $r$, $s$, bringing the product
$$\begin{pmatrix}
1 & 0 & 0 \\
q & 1 & 0 \\
r & s & 1 \end{pmatrix}\cdot \begin{pmatrix}
1 & 0 & 0 \\
0 & t^b & 0 \\
0 & 0 & t^c
\end{pmatrix}$$
into one of the forms given in \S\ref{markergerms}.
Thus, once coordinates are chosen in the target plane,
contributing germs are
equivalent to germs of the form
$$\alpha_\bullet(t) \cdot M$$
where $\alpha_\bullet(t)=\alpha_{III}(t)$, $\alpha_{IV}(t)$, or $\alpha_V(t)$,
according to the type of the component. Replacing $t$ by a root $t^{1/k}$
ensures, if necessary, that the exponents appearing in such a germ are 
relatively prime.
\end{remark}

\subsection{Computation of the multiplicities}\label{multin}

\begin{defin}
The {\em weight\/} of a germ $\alpha(t)$ in $\Pbb^8$ is the order of 
vanishing in $t$ of $F\circ\alpha(t)$.
\end{defin}

Note that the weight of $\alpha(t)$ is the order of contact of
$\alpha(t)$ with $\cS$: indeed, it is the minimum intersection
multiplicity of $\alpha(t)$ and generators
$F\circ\varphi(x_0,y_0,z_0)$ of the ideal of~$\cS$, at $\alpha(0)$
(cf.~\S\ref{beginprel}).

\begin{lemma}\label{weightsmult}
The multiplicity $m_{ij}$ is the minimum weight of a germ $\alpha(t)$
such that $\oalpha\in \oE_{ij}$.
\end{lemma}

\begin{proof}
Let $\oalpha$ be a general point of $\oE_{ij}$. Since $\nG$ is normal,
we may assume that it is nonsingular at $\oalpha$. Let $(\overline
z)=(z_1,\dots,z_8)$ be a system of local parameters for $\nG$ centered
at $\oalpha$, and such that the ideal of $\oE_{ij}$ is $(z_1)$ near
$\oalpha$; thus the ideal of $\oE$ is $(z_1^{m_{ij}})$ near~$\oalpha$. 
Consider the germ $\oalpha(t)$ in $\nG$ defined by
$$\oalpha(t)=(t,0,\dots,0)\quad,$$
and its push-forward $\alpha(t)=\onu(\oalpha(t))$ in $\Pbb^8$.

The weight of $\alpha(t)$ is the order of contact of $\alpha(t)$ with
$\cS$; hence it equals the order of contact of $\oalpha(t)$ with
$\onu^{-1}(\cS)=\oE$; pulling back the ideal of $\oE$ to $\oalpha(t)$,
we see that this equals $m_{ij}$. 

Any other germ in $\nG$ meeting $\oE_{ij}$ and not contained in $\oE$
must have intersection number $\ge m_{ij}$; the statement follows.
\end{proof}

In fact, a germ in $\Pbb^8$ that lifts to a germ in $\nG$ meeting the 
support of $\oE$ at a general point of $\oE_{ij}$ intersects $\oE_{ij}$
{\em transversally\/} if and only if its weight is $m_{ij}$.

Now, by Lemma~3.2 in \cite{ghizzI},
equivalent germs $\alpha(t)$, $\beta(t)$ have the same weight; and their 
lifts are centered at the same point $\oalpha=\obeta$, by  
Lemma~\ref{equivgermso}. 
In particular, if $\alpha(t)$ and $\beta(t)$ are equivalent, then 
$\alpha(t)$ lifts to a germ transversal to a component $\oE_{ij}$ at a 
general point if and only if $\beta(t)$ does.

The following easy consequence of these considerations yields the 
multiplicities $m_{ij}$:

\begin{corol}\label{weightcor}
The multiplicity $m_{ij}$ of a component $\oE_{ij}$ of $\oE$ over a 
component $E_i$ equals the weight of a corresponding {\em marker\/} 
germ.
\end{corol}

\begin{proof}
By Lemma~\ref{weightsmult}, there exists a germ $\beta(t)$ of
weight $m_{ij}$, meeting $\oE_{ij}$ at a (general) point $\obeta$.
Such a germ $\beta(t)$ contributes the component $E_i$, hence
$\beta(t)$ is equivalent to $\alpha(t^k)$ for a corresponding marker 
germ $\alpha(t)$, by Lemma~\ref{markequiv}. The weight of $\alpha(t)$
is then $m_{ij}/k$, by Lemma~3.2 in \cite{ghizzI}. Since $m_{ij}$
is the minimum weight, we have $k=1$, and the weight of
the marker germ $\alpha(t)$ is $m_{ij}$ as stated.
\end{proof}

It is now straightforward to verify the multiplicity statements given in 
Propositions~\ref{typeIIImul}---\ref{typeVmul}. For 
example, for type~IV: with notation as above (and in \cite{ghizzI}, \S2.3,
with $q=e$, ${\overline q}={\overline e}$, $r=f$),
the initial term of
$\cC\circ \alpha_{IV}(t)$ is
$$x^{\overline q}y^rz^q \prod_{j=1}^S(y^c+\rho_j x^{c-b}z^b)\,
t^{Sbc+br+cq}\quad;$$
as $(j_0,k_0)=(r,q+Sb)$, $(j_1,k_1)=(r+Sc,q)$, the weight is
$$Sbc+br+cq=\frac{j_1k_0-j_0k_1}S$$
as stated.
The other two verifications are left to the reader; for type~V, use
Lemma~4.3 in~\cite{ghizzI}.

\section{Components in the normalization}\label{compnor}
\subsection{The \textbf{PGL(3)}-action}\label{pgl3act}
The other information listed in Propositions~\ref{typeIIImul},
\ref{typeIVmul}, and~\ref{typeVmul} requires a more explicit
description of the components $\oD$ of $\oE$, especially in
connection with the behavior of germs centered on these components.

One important ingredient is the $\PGL(3)$-action on $\nG$.
The $\PGL(3)$-action on $\Pbb^8$ given by multiplication on the right
makes the basic rational map
$\Pbb^8 \dashrightarrow \Pbb^N$
equivariant, and hence induces a right $\PGL(3)$-action on
$\TPbb^8$ and $\nG$, fixing each component of~$\oE$. Explicitly,
if $\oalpha$ is the center of the lift to $\nG$ of a germ $\alpha(t)$,
then $\oalpha\cdot N$ is defined to be the center of the lift of
the germ $\alpha(t)\cdot N$, for $N\in \PGL(3)$. We record the
following trivial but useful remarks:

\begin{lemma}\label{actionlift}
If $\beta(t)=\alpha(t)\cdot N$ for $N\in\PGL(3)$,
and $\oalpha$ belongs to a component $\oD$ of $\oE$, then
$\obeta\in \oD$.
\end{lemma}

\begin{lemma}\label{formal}
Let $\oD$ be a component of $\oE$ over a component $D$ of type 
III,~IV, or~V.  Then the orbit of a general $\oalpha\in\oD$ is dense 
in $\oD$.
\end{lemma}

\noindent (Indeed, the orbit of a general point in a component of type
III,~IV, or~V has dimension~7, as follows from the explicit description
of such points given in \S\ref{prell}.)

Henceforth, $D$ will denote a component of $E$ of type~III, IV, or~V;
and $\oD$ will be a component of $\oE$ over $D$. A {\em general 
point\/} of $\oD$ (resp.~$D$) will be a point of the dense $\PGL(3)$ 
orbit in $\oD$ (resp.~$D$). A germ $\alpha(t)$ `contributes' to $D$ if 
$\oalpha$ is a general point of $\oD$ in this sense.

\subsection{Criterion for equal lift}\label{cfel}
We are ready to upgrade Lemma~\ref{equivgermso}.

\begin{prop}\label{criteqlift}
Let $\alpha(t)$, $\beta(t)$ be germs such that $\oalpha(t)$, $\obeta(t)$ 
are centered at general points of components of $\oE$ dominating
components of type~III, IV, or~V, and meet these components transversally. 
Then $\oalpha=\obeta$ if and only if $\alpha(t)$, $\beta(t)$ are equivalent.
\end{prop}

\begin{proof}
One implication is given in Lemma~\ref{equivgermso}.

We will give the argument for the converse under the assumption
that the entries in $\alpha(t)$ are polynomials; this
is the only case in which we will use the statement, and we leave
to the reader the (easy) extension to the general case.

Assume $\oalpha=\obeta$ is a general point of a component $\oD$
of $\oE$, and $\oalpha(t)$, $\obeta(t)$ meet $\oD$ transversally.

The image of $\oalpha=\obeta$ in $\TPbb^8$ is a point $(\alpha(0),
\cX)$ of a unique component $D$ of $E$ of type III,~IV, or~V;
note that the $\PGL(3)$-stabilizer of $(\alpha(0),\cX)$ has dimension~1. 
Consider an $\Abb^7\subset
\Pbb^8$ through the identity $I$ and transversal at $I$ to the stabilizer;
let $U=\Abb^7\cap\PGL(3)$, and consider the action map 
$\Abb^1\times U\to \nG$: 
$$(t,\varphi) \mapsto \oalpha(t)\circ\varphi\quad.$$
This map is dominant, and \'etale at $(0,I)$.
Note that $\alpha(t)$ factors through it:
$$\xymatrix@R=5pt@C=15pt{
\Abb^1 \ar[r] & \Abb^1\times U \ar[r] & \nG \ar[r] & \Pbb^8 &\\
t \ar@{|->}[r] & (t,I) \ar@{|->}[r] & \oalpha(t) \ar@{|->}[r] & \alpha(t) &.
}$$
Parametrizing a lift of $\beta(t)$ to $\Abb^1\times U$ we likewise get a 
factorization
$$t \mapsto (z(t),M(t)) \mapsto \oalpha(z(t))\circ M(t)=\obeta(t) 
\mapsto \beta(t)$$
for suitable ($\Cbb[[t]]$-valued) $z(t)$, $M(t)$.
Since $\oalpha=\obeta$ in $\nG$,
we may assume that the center $(z(0),M(0))$ of the lift of $\beta(t)$ 
equals the center $(0,I)$ of the lift of $\alpha(t)$. Also, $z(t)$ 
vanishes to order~1 at $t=0$, since $\obeta(t)$
is transversal to $\oD$.
Hence there exists a unit $\nu(t)$ such that $z(t\nu(t))=t$, and we
can apply the parameter change
$$\beta(t\nu(t))=\alpha(t)\circ M(t\nu(t))=\alpha(t)\circ m(t)\quad,$$
where we have set $m(t)=M(t\nu(t))$, a $\Cbb[[t]]$-valued point of 
$\PGL(3)$.

As $m(0)=M(0)=I$, this shows that $\alpha(t)$ and $\beta(t)$ are
equivalent in the sense of Definition~\ref{equivgermsnew}, concluding the 
proof.
\end{proof}

\begin{corol}\label{critesamecom}
Let $\alpha(t)$, $\beta(t)$ be germs such that $\oalpha(t)$, $\obeta(t)$ 
are centered at general points of components of $\oE$ dominating
components of type~III, IV, or~V, and meet these components transversally. 

Then $\oalpha$, $\obeta$ belong to the same component of $\oE$ if
and only if $\alpha(t)^{-1}\beta(\tau(t))$ is a $\Cbb[[t]]$-valued point of
$\PGL(3)$, for a change of parameter $\tau(t)=t\nu(t)$ with
$\nu(t)\in\Cbb[[t]]$ a unit.
\end{corol}

\begin{proof}
This is an immediate consequence of Proposition~\ref{criteqlift} and of 
Lemma~\ref{formal}.
\end{proof}

\subsection{Technical lemma}
Applications of Proposition~\ref{criteqlift} and Corollary~\ref{critesamecom} 
to the case of type~V components will rely on the following technical lemma.

\begin{lemma}\label{technical}
Let 
$$\alpha(t)=\begin{pmatrix}
1 & 0 & 0 \\
t^a & t^b & 0 \\
\underline{f(t^a)} & \underline{f'(t^a)t^b} & t^c 
\end{pmatrix}\quad,\quad
\beta(t)=\begin{pmatrix}
1 & 0 & 0 \\
t^{a'} & t^{b'} & 0 \\
\underline{g(t^{a'})} & \underline{g'(t^{a'})t^{b'}} & t^{c'} 
\end{pmatrix}$$
be two marker germs for type~V components, and assume that
$\alpha(t)^{-1}\beta(\tau(t))$ is a $\Cbb[[t]]$-valued point of
$\PGL(3)$, for a change of parameter $\tau(t)=t\nu(t)$ with
$\nu(t)\in\Cbb[[t]]$ a unit. Then $a'=a$, $b'=b$, $c'=c$, and
$\nu(t)=\xi(1+t^{b-a}\mu(t))$, where $\xi$ is an $a$-th root of~$1$
and $\mu(t)\in\Cbb[[t]]$; further, $\underline{g((\xi
t)^a)}=\underline{f(t^a)}$.
\end{lemma}

Recall that $\underline{\dots}$ stands for a truncation; modulo $t^c$
in the first germ, and modulo~$t^{c'}$ in the second germ. The notation
is unambiguous once the equality $(a,b,c)=(a',b',c')$ is established.
Also, the expression $g((\xi t)^a)$ is then an abbreviation for 
$\sum \xi^{a\lambda_i} \gamma'_{\lambda_i} t^{a\lambda_i}$, where 
$g(y)=\sum \gamma'_{\lambda_i} y^{\lambda_i}$; the coefficients 
$\xi^{a\lambda_i}$ are well-defined for $\lambda_i<C$ since 
$a\lambda_i$ is an integer for $\lambda_i<C$ (cf.~\S\ref{result}).

\begin{proof}
Write $\varphi(t)=\underline{f(t^a)}$ and
$\psi(t)t^b=\underline{f'(t^a)t^b}$. The hypothesis is that
$$\alpha(t)^{-1}\cdot\beta(\tau)=\begin{pmatrix}
1 & 0 & 0 \\
\frac{\tau^{a'}-t^a}{t^b} & \frac{\tau^{b'}}{t^b} & 0 \\
\frac{\underline{g(\tau^{a'})}-\varphi(t)-(\tau^{a'}-t^a)\psi(t)}
     {t^c} &
\frac{\underline{g'(\tau^{a'})\tau^{b'}}-\psi(t)\tau^{b'}}{t^c} &
\frac{\tau^{c'}}{t^c}
\end{pmatrix}$$
has entries in $\Cbb[[t]]$, and its determinant is a unit in
$\Cbb[[t]]$. This implies $b'=b$ and $c'=c$. As 
$$\frac{\tau^{a'}-t^a}{t^b}\in \Cbb[[t]]$$
and $b>a$, necessarily $a'=a$ and $t^a(\nu(t)^a-1)=(\tau^a-t^a)\equiv
0\mod{t^b}$. This implies 
$$\nu(t)=\xi(1+t^{b-a}\mu(t))$$ 
for $\xi$ an $a$-th root of $1$ and $\mu(t)\in\Cbb[[t]]$. 
Also note that since the triples $(a,b,c)$ and $(a',b',c')$ coincide,
necessarily the dominant term in $g(y)$ has the same exponent
$\lambda_0$ as in $f(y)$,
since $a\lambda_0=2a-2b+c$ (see \S\ref{markergerms}).

Now we claim that
$$g(\tau^a)-(\tau^a-t^a)\psi(t)\equiv g((\xi t)^a)\mod t^c\quad.$$
Granting this for a moment, it follows that
$$\underline{g(\tau^{a'})}-\varphi(t)-(\tau^{a'}-t^a)\psi(t) \equiv
g((\xi t)^a)-f(t^a)\mod t^c\quad;$$
hence, the fact that the $(3,1)$ entry is in $\Cbb[[t]]$ implies that
$$g((\xi t)^a)\equiv f(t^a)\mod t^c\quad,$$
which is what we need to show in order to complete the proof.

Since the $(3,2)$ entry is in $\Cbb[[t]]$, necessarily
$$g'(\tau^a)\equiv \psi(t)\mod {t^{c-b}}\quad;$$
so our claim is equivalent to the assertion that
$$g(\tau^a)-(\tau^a-t^a)g'(\tau^a)\equiv g((\xi t)^a)\mod t^c\quad.$$
By linearity, in order to prove this it is enough to verify the
stated congruence for $g(y)=y^\lambda$, with $\lambda\ge
\lambda_0$. That is, we have to verify that if $\lambda\ge \lambda_0$
then
$$\tau^{a\lambda}-(\tau^a-t^a)\lambda \tau^{a\lambda-a}\equiv (\xi
t)^{a\lambda} \mod t^c\quad.$$
For this, observe
$$\tau^{a\lambda}=(\xi t)^{a\lambda}(1+t^{b-a}\mu(t))^{a\lambda}
\equiv(\xi t)^{a\lambda}(1+a\lambda t^{b-a}\mu(t)) \mod
      {t^{a\lambda+2(b-a)}}$$
and similarly
$$\tau^{a\lambda-a}=(\xi t)^{a\lambda-a}(1+t^{b-a}\mu(t))^{a\lambda-a}
\equiv t^{-a}(\xi t)^{a\lambda}\mod {t^{a\lambda-a+(b-a)}}\quad,$$
$$(\tau^a-t^a)=(\xi t)^a(1+t^{b-a}\mu(t))^a-t^a\equiv at^b\mu(t)
\mod {t^{a+2(b-a)}}\quad.$$
Thus
$$(\tau^a-t^a)\lambda \tau^{a\lambda-a}\equiv (\xi t)^{a\lambda}
a\lambda t^{b-a}\mu(t)\mod {t^{a\lambda+2(b-a)}}$$
and
$$\tau^{a\lambda}-(\tau^a-t^a)\lambda \tau^{a\lambda-a}\equiv (\xi
t)^{a\lambda} \mod {t^{a\lambda+2(b-a)}}\quad.$$
Since
$$a\lambda+2(b-a)\ge a\lambda_0+2b-2a=c\quad,$$
our claim follows.
\end{proof}

\subsection{Number of components in the normalization}
We are ready to prove the statements in Propositions~\ref{typeIIImul},
\ref{typeIVmul}, and~\ref{typeVmul} concerning the number of
components $\oD$ over a given component $D$ of $E$ of type~III,
IV, or~V.

\textbf{Type~III and~IV.} For type~III and~IV components, the statement
is that for any fixed $p\in \cC$ and (for type~IV) line $L$ in the tangent
cone to $p$ at $\cC$, there is exactly one component $\oD$ over each 
component $D$ of $E$.

This is in fact an easy consequence of Lemma~\ref{actionlift}.
For example, in the case of components of type~IV it suffices to verify that, 
for fixed $L$,
any two marker germs for a given component
lift to germs in $\nG$ centered on the same component of $\oE$. 
Now, such marker germs are of the form
(cf.~\S\ref{markergerms} and Remark~\ref{coordchoice})
$$
\begin{pmatrix}
1 & 0 & 0 \\
0 & t^b & 0 \\
0 & 0 & t^c
\end{pmatrix}\cdot M_1\quad,\quad
\begin{pmatrix}
1 & 0 & 0 \\
0 & t^b & 0 \\
0 & 0 & t^c
\end{pmatrix}\cdot M_2
$$
for two invertible matrices $M_1$, $M_2$. Lemma~\ref{actionlift} implies
immediately that the lifts of these two germs are centered on the same
component of $\oE$.\qed

\textbf{Type~V.} The situation for components of type~V is more complex,
and requires the use of Corollary~\ref{critesamecom} and 
Lemma~\ref{technical}.

For a fixed point $p$ and line $L$, and once coordinates are chosen
as usual
(so that $p=(1:0:0)$, and $L$ is the line $z=0$) type~V components are 
determined by pairs $(C,f(y))$, where $z=f(y)$ is a formal branch of
$\cC$, by the procedure described in \S\ref{result}. 
Recall that $(C,f(y))$ and $(C',g(y))$ are {\em sibling\/} 
data if  the corresponding integers $a<b<c$, $a'<b'<c'$ are the same (so in
particular $C=C'$) and further 
$$g_{(C)}(y)=\sum_{\lambda_i<C} \xi^{a\lambda_i} \gamma_{\lambda_i}
y^{\lambda_i}$$
for an $a$-th root $\xi$ of $1$ (or, in abbreviated form, $g_{(C)}(t^a)=
f_{(C)}((\xi t)^a)$).
The statement we must prove is the following:

\begin{claim}\label{siblingsame}
Two pairs $(C,f(y))$ and $(C',g(y))$ determine the same
component $\oD$ over $D$ if and only if they are siblings.
\end{claim}

\begin{proof}
Let $\alpha(t)$, $\beta(t)$ be two marker germs leading to $D$;
we may assume (cf.~Remark~\ref{coordchoice}, 
Lemma~\ref{actionlift}) that
$$\alpha(t)=\begin{pmatrix}
1 & 0 & 0 \\
t^a & t^b & 0 \\
\underline{f(t^a)} & \underline{f'(t^a)t^b} & t^c 
\end{pmatrix}
\quad,\quad
\beta(t)=\begin{pmatrix}
1 & 0 & 0 \\
t^{a'} & t^{b'} & 0 \\
\underline{g(t^{a'})} & \underline{g'(t^{a'})t^{b'}} & t^{c'} 
\end{pmatrix}\quad.$$
If these two germs determine the same component of $\oE$, then 
by Corollary~\ref{critesamecom} $\alpha(t)^{-1}\beta(t\nu(t))$ is 
a $\Cbb[[t]]$-valued point of $\PGL(3)$, for a unit $\nu(t)\in\Cbb[[t]]$.
It follows (by Lemma~\ref{technical}) that $a'=a$, $b'=b$, $c'=c$, and
$\underline{g((\xi t)^a)}=\underline{f(t^a)}$; that is, $(C,f(y))$ and
$(C',g(y))$ are siblings.

Conversely, assume that $(C,f(y))$, $(C',g(y))$ are siblings. Then 
$C=C'$, and for an $a$-th root $\xi$ of $1$ the corresponding germs
$$\alpha(t)=\begin{pmatrix}
1 & 0 & 0 \\
t^a & t^b & 0 \\
\underline{f(t^a)} & \underline{f'(t^a)t^b} & t^c 
\end{pmatrix}
\quad,\quad
\beta(t)=\begin{pmatrix}
1 & 0 & 0 \\
t^a & t^b & 0 \\
\underline{g(t^a)} & \underline{g'(t^a)t^b} & t^c 
\end{pmatrix}$$
satisfy 
$$\aligned
\beta(\xi t) &=\begin{pmatrix}
1 & 0 & 0 \\
(\xi t)^a & (\xi t)^b & 0 \\
\underline{g((\xi t)^a)} & \underline{g'((\xi t)^a)(\xi t)^b} & 
(\xi t)^c 
\end{pmatrix} = \begin{pmatrix}
1 & 0 & 0 \\
t^a & t^b\xi^b  & 0 \\
\underline{f(t^a)} & \underline{f'(t^a)t^b}\xi^b  & t^c \xi^c 
\end{pmatrix}\\
& = \alpha(t) \cdot \begin{pmatrix}
1 & 0 & 0 \\
0 & \xi^b & 0 \\
0 & 0 & \xi^c
\end{pmatrix}\quad.
\endaligned$$
Therefore, $\alpha(t)^{-1}\beta(\xi t)$ is an invertible constant matrix.
This shows that $\oalpha$ and $\obeta$ belong to the same component
of $\oE$, by Lemma~\ref{actionlift}.
\end{proof}

This concludes the proof of the statement concerning the number of
components~$\oD$ over a given component $D$ of $E$ in 
Propositions~\ref{typeIIImul}, \ref{typeIVmul}, \ref{typeVmul}.

\section{End of the proof}\label{Eop}
\subsection{Inessential subgroups}
We are left with the task of verifying the statement concerning the
degrees of the maps $\oD\to D$; our main tool here will be 
Proposition~\ref{criteqlift}.

Recall that $\PGL(3)$ acts on both $\oD$ and on the underlying component
$D$. Accordingly, every general $\oalpha\in \oD$ determines two 
one-dimensional subgroups of $\PGL(3)$:
\begin{itemize}
\item the $\PGL(3)$-stabilizer $\Stab\oalpha$ of $\oalpha$; and
\item the $\PGL(3)$-stabilizer $\Stab{(\alpha,\cX)}$ of the image of
$\oalpha$ in $D$.
\end{itemize}
The equivariance of the normalization map $n:\nG \to \TPbb^8$ implies that
$\Stab\oalpha$ is a subgroup of (in fact, a union of components of) 
$\Stab{(\alpha,\cX)}$.

\begin{lemma}\label{stablemma}
Let $D$ be a component of type~III, IV, or~V, and let
$\alpha(t)$ be a marker germ for $D$. Let $\alpha=\alpha(0)$, 
$\cX=\lim \cC\circ \alpha(t)$, and let $\oD$ be the component
of $\oE$ over~$D$ containing the center $\oalpha$ of the lift of
$\alpha(t)$ to $\nG$. 

Then the degree of $\oD$ over $D$ is the index of $\Stab\oalpha$ in 
$\Stab{(\alpha,\cX)}$.
\end{lemma}

\begin{proof}
If $\alpha(t)$ is a marker germ then $(\alpha,\cX)$ and $\oalpha$ are
general; by Lemma~\ref{formal}, $\oD$ is the closure of the 
$\PGL(3)$-orbit of $\oalpha$. It follows that the fiber of $\oD\to D$
over $(\alpha,\cX)$ is the $\Stab{(\alpha,\cX)}$-orbit of $\oalpha$,
giving the statement.
\end{proof} 

The stabilizers $\Stab{(\alpha,\cX)}$ are easily identified subgroups of
the stabilizers of the curves $\cX$, which are discussed in 
\cite{MR2002d:14084}, \S1. We have to determine the stabilizers
$\Stab\oalpha$, and we do this by means of the following construction.

Let $\alpha(t)$ be a marker germ whose lift to $\nG$ is centered at
$\oalpha\in \oD$. Consider the $\Cbb((t))$-valued points of $\PGL(3)$ 
obtained as products 
$$M_\nu(t):=\alpha(t)^{-1}\cdot \alpha(t \nu(t) )$$
as $\nu(t)$ ranges over all units in $\Cbb[[t]]$. 
Among all the $M_\nu(t)$, consider those that are in fact 
$\Cbb[[t]]$-valued points of $\PGL(3)$, and in that case let
$$M_\nu:=M_\nu(0)\quad.$$

The reader may verify directly that the set of all $M_\nu$ so obtained is a 
subgroup of the stabilizer of $(\alpha,\cX)$. We call it the {\em inessential
subgroup\/} (w.r.t.~$\alpha(t)$) of the stabilizer of $(\alpha,\cX)$.
It consists of the elements of the stabilizer due to  
reparametrizations of the germ $\alpha(t)$.
We now prove that this subgroup equals the stabilizer of $\oalpha$.

\begin{prop}\label{degreetool}
Let $D$ be a component of $E$ of type~III, IV, or~V, and let $\oD$ 
be any component of $\oE$ dominating $D$. Further, let $\alpha(t)$
be a marker germ for $D$, such that the lift of $\alpha(t)$ to $\nG$
is centered at a point $\oalpha\in \oD$. Let $\alpha=\alpha(0)$, 
$\cX=\lim \cC\circ \alpha(t)$.

Then the inessential subgroup of $\Stab{(\alpha,\cX)}$ 
{\em ({\em w.r.t.~$\alpha(t)$})} is the stabilizer
of $\oalpha$, and the degree of $\oD$ over $D$ is the index of the
inessential subgroup in $\Stab{(\alpha,\cX)}$.
\end{prop}

\begin{proof}
Let $M_\nu$ (as above) be an element of the inessential subgroup. Then
$\alpha(t) \cdot M_\nu$ and $\alpha(t)\cdot M_\nu(t)=\alpha(t\nu(t))$ are 
equivalent according to Definition~\ref{equivgermsnew}; by 
Lemma~\ref{equivgermso}, $\oalpha\cdot M_\nu=\oalpha$. 

For the converse, assume $\oalpha=\oalpha\cdot N$. By
Proposition~\ref{criteqlift}, $\alpha(t)$ is equivalent to $\alpha(t)\cdot 
N$: that is, there is a $\Cbb[[t]]$-valued point $N(t)$ of $\PGL(3)$,
with $N(0)=N$, and a unit $\nu(t)\in \Cbb[[t]]$, such that 
$$\alpha(t\nu(t))= \alpha(t)\cdot N(t)\quad.$$
Therefore $M_\nu(t)=\alpha(t)^{-1}\cdot \alpha(t\nu(t))=N(t)\in \Cbb[[t]]$: 
that is, $N$ is in the inessential subgroup of the stabilizer of $(\alpha,\cX)$.

The statement about the degree of $\oD$ over $D$ follows from
Lemma~\ref{stablemma}, completing the proof.
\end{proof}

\subsection{The degree of $\oD$ over $D$}\label{degoDoD}
We are ready to complete the proof of Propositions~\ref{typeIIImul},
\ref{typeIVmul}, and~\ref{typeVmul}, and hence of Theorem~\ref{main}.

All that is left to prove is the statement concerning the degree of each
component of $\oE$ over the corresponding component of $E$; this
is done by repeated applications of Proposition~\ref{degreetool},
that is, by determining the inessential subgroups of the
stabilizers for components of type~III, IV, and~V.

\begin{claim}
For type~III and~IV, the inessential subgroup is the component of
the identity in the stabilizer of $(\alpha,\cX)$.
\end{claim}

\begin{proof}
For type~III and~IV, marker germs are of the form
$$\alpha(t)=\begin{pmatrix}
1 & 0 & 0 \\
0 & t^b & 0 \\
0 & 0 & t^c
\end{pmatrix}\quad,$$
with $b=c=1$ for type~III and $b$, $c$ positive and relatively prime
for type~IV.

For all units $\nu(t)$,
$$\alpha(t)^{-1}\cdot \alpha(t\nu(t))=
\begin{pmatrix}
1 & 0 & 0 \\
0 & t^b & 0 \\
0 & 0 & t^c
\end{pmatrix}^{-1}\cdot\begin{pmatrix}
1 & 0 & 0 \\
0 & t^b\nu(t)^b & 0 \\
0 & 0 & t^c\nu(t)^c
\end{pmatrix}=
\begin{pmatrix}
1 & 0 & 0 \\
0 & \nu(t)^b & 0\\
0 & 0 & \nu(t)^c
\end{pmatrix}$$
is a $\Cbb[[t]]$-valued point of $\PGL(3)$. Thus the inessential
subgroup consists of the elements
$$\begin{pmatrix}
1 & 0 & 0 \\
0 & \nu^b & 0\\
0 & 0 & \nu^c
\end{pmatrix}\quad,$$
for $\nu=\nu(0)\in\Cbb$, $\nu\ne 0$. These elements form 
the identity component of the stabilizer (cf.~\S3
in \cite{MR2002d:14083}), verifying our claim.
\end{proof}

Applying Proposition~\ref{degreetool}, we conclude that for
types~III,~IV the degree of $\oD\to D$ equals the number
of components of the stabilizer of a general point $(\alpha,\cX)
\in D$. 

\textbf{Type~III.} Recall that the limit of $\cC$ along a marker germ 
$\alpha(t)$ consists of a fan $\cX$ whose star reproduces 
the tangent cone to $\cC$ at $p$, and whose free line is supported on 
the kernel line $x=0$. It is easily checked that the stabilizer of
$(\alpha(0),\cX)$ has one component for each element
of $\PGL(2)$ fixing the $m$-tuple determined by the tangent cone to
$\cC$ at $p$, verifying the degree statement in 
Proposition~\ref{typeIIImul}.\qed

\textbf{Type~IV.} The number of components of the stabilizer of a general
$(\alpha,\cX)$ is determined as follows. The limit $\cX$ is given by
$$x^{\overline q}y^rz^q \prod_{j=1}^S(y^c+\rho_j x^{c-b}z^b)\quad;$$
the stabilizer of $(\alpha,\cX)$ is the subgroup of the stabilizer of $\cX$ 
fixing the kernel line $x=0$. 
Thus, the number of components of the stabilizer of $(\alpha,\cX)$
equals the number of components of the stabilizer of $\cX$, or
the same number divided by~2, according to whether the kernel line is
determined by $\cX$ or not; the latter eventuality occurs
precisely when $c=2$ and $q=\overline q$. It follows then from 
Lemma~3.1 in \cite{MR2002d:14083} that the number of components
of the stabilizer of $(\alpha,\cX)$ equals the number of 
automorphisms $\Abb^1 \to \Abb^1$, $\rho\mapsto u\rho$ (with $u$ a 
root of unity) preserving the $S$-tuple
 $[\rho_1,\dots,\rho_S]$.
This completes the proof of Proposition~\ref{typeIVmul}.\qed
\vskip 6pt

\textbf{Type~V.}
Finally, we deal with components of type~V. The determination of
the inessential subgroup of the stabilizer of a general point
$(\alpha,\cX)$ of such a component relies again on the
technical Lemma~\ref{technical}.

\begin{lemma}\label{typeVines}
Let
$$\alpha(t)=\begin{pmatrix}
1 & 0 & 0 \\
t^a & t^b & 0 \\
\underline{f(t^a)} & \underline{f'(t^a)t^b} & t^c 
\end{pmatrix}$$
be the marker germ determined by $C$ and the formal branch
$f(y)=\sum \gamma_{\lambda_i}y^{\lambda_i}$, and let
$\cX= \lim_{t\to 0} C\circ\alpha(t)$.
Then the corresponding inessential subgroup 
consists of the components of the stabilizer of 
$(\alpha,\cX)$ containing 
matrices
$$\begin{pmatrix}
1 & 0 & 0 \\
0 & \eta^b & 0\\
0 & 0 & \eta^c
\end{pmatrix}$$
with $\eta$ an $h$-th root of~$1$, where $h$ is the greatest common
divisor of $a$ and all $a\lambda_i$ $(\lambda_i<C)$.
\end{lemma}

\begin{proof}
For every $h$-th root $\eta$ of~$1$, each component of the stabilizer
containing a diagonal matrix of the given form is in the inessential
subgroup: indeed, such a diagonal matrix can be realized as
$\alpha(t)^{-1}\cdot \alpha(\eta t)$.

To see that, conversely, every component of the inessential subgroup
is as stated, apply Lemma~\ref{technical} with
$\beta(t)=\alpha(t)$. We find that if $\alpha(t)^{-1}\cdot
\alpha(t\nu(t))$ is a $\Cbb[[t]]$-valued point of $\PGL(3)$, then
$\nu(t)=\eta(1+t^{b-a}\mu(t))$, with $\eta$ an $a$-th root of~$1$, and
further 
$$\underline{f(t^a)}=\underline{f((\eta t)^a)}\quad,$$
that is,
$$\sum_{\lambda_i<C} \gamma_{\lambda_i} y^{\lambda_i}
=\sum_{\lambda_i<C} \eta^{a\lambda_i}\gamma_{\lambda_i}
y^{\lambda_i}\quad.$$ 
Therefore $\eta^{a\lambda_i}=1$ for all $i$ such that $\lambda_i<C$,
and it follows that $\eta$ is an $h$-th root of~$1$. 

For $\nu(t)=\eta(1+t^{b-a}\mu(t))$, the matrix $M_\nu(0)=\alpha(t)^{-1}\cdot
\alpha(t\nu(t))|_{t=0}$ is lower triangular and invertible, of the form 
$$\begin{pmatrix}
1 & 0 & 0 \\
a\mu_0 & \eta^b & 0\\
 \gamma_{\lambda_0}\binom{\lambda_0}2 (a\mu_0)^2 
+ \gamma_{\frac{\lambda_0+C}2} \frac{\lambda_0+C}2 (a\mu_0) &
2\gamma_{\lambda_0}\binom{\lambda_0}2 (a\mu_0) \eta^b & \eta^c
\end{pmatrix}$$
where $\mu_0=\mu(0)$. These matrices are in the stabilizer of
$(\alpha,\cX)$ for all $\mu_0$ (since they are in the inessential subgroup).
Setting $\mu_0=0$ proves the statement.
\end{proof}

Note that $\eta^c=(\eta^b)^2$ since $c-2b=a\lambda_0-2a$ is divisible by $h$;
this is in fact a necessary condition for the diagonal matrix above to belong
to the stabilizer. Moreover, if $\gamma_{\frac{\lambda_0+C}2}\neq0$, then
necessarily $\eta^b=1$; as the proof of the following proposition shows,
this implies $h=1$.

\begin{prop}\label{typeVmuldeg}
For the component $\oD$ determined by the truncation $f_{(C)}(y)$ as
above, 
let $A$ be the number of components of the stabilizer of the limit
$$x^{d-2S}\prod_{i=1}^S\left(zx-\frac {\lambda_0(\lambda_0-1)}2
\gamma_{\lambda_0}y^2 -\frac{\lambda_0+C}2
\gamma_{\frac{\lambda_0+C}2}yx-\gamma_C^{(i)}x^2\right)$$
(that is, by \cite{MR2002d:14083}, \S4.1, twice the number of
automorphisms $\gamma\to u\gamma+v$ preserving the $S$-tuple
$[\gamma_C^{(1)},\dots,\gamma_C^{(S)}]$).
Then the degree of the map $\oD \to D$ equals $\frac Ah$, where $h$ is
the number determined in Lemma~\ref{typeVines}.
\end{prop}

\begin{proof}
As the kernel line must be supported on the distinguished
tangent of the limit~$\cX$, the stabilizer of $(\alpha,
\cX)$ equals the stabilizer of $\cX$, and in particular
it consists of $A$ components.

Next, observe that for $\eta_1\ne \eta_2$ two $h$-th roots of~$1$, the
two matrices
$$\begin{pmatrix}
1 & 0 & 0 \\
0 & \eta_1^b & 0\\
0 & 0 & \eta_1^c
\end{pmatrix}\quad,\quad
\begin{pmatrix}
1 & 0 & 0 \\
0 & \eta_2^b & 0\\
0 & 0 & \eta_2^c
\end{pmatrix}$$
are distinct: indeed, if $\eta^b=\eta^c=1$, then the order of $\eta$
divides every exponent of every entry of $\alpha(t)$, hence it equals~$1$
by the minimality of~$a$. Further, the components of the stabilizer
containing these two matrices must be distinct: indeed, the
description of the identity component of the stabilizer of a curve
consisting of quadritangent conics given in \cite{MR2002d:14084}, \S1,
shows that the only diagonal matrix in the component of the identity
is in fact the identity itself. 

Therefore the index of the inessential subgroup equals $A/h$, and the
statement follows then from Proposition~\ref{degreetool}.
\end{proof}

Proposition~\ref{typeVmuldeg} verifies the degree statement in
Proposition~\ref{typeVmul}, thereby completing the proof of that proposition,
and hence of Theorem~\ref{main}.

\bibliographystyle{abbrv}
\bibliography{ghizzIIbib}

\end{document}